\documentclass[11pt]{article}
\usepackage{latexsym,amssymb,amsmath,amscd,amsthm,amsxtra}

\newcommand{\counte}{section}

\author{Zhao Xu-an, zhaoxa@bnu.edu.cn\\ Department of Mathematics, Beijing Normal University\\Key Laboratory
of Mathematics and Complex Systems\\ Ministry of Education,
China, Beijing 100875\\Jin Chunhua, jinch@amss.ac.cn\\Academy of Mathematics and Systems Science, Chinese Academy of Sciences
\\China, Beijing 100190
\\ Zhang Jimin\\Department of Mathematics, Beijing Normal University}

\title{Poincar\'{e} series and rational cohomology rings of Kac-Moody groups and their flag manifolds\thanks{The authors are supported by NSFC11171025}}

\begin{document}

\maketitle
\begin{abstract}
In this paper, we study the rational cohomology rings of indefinite Kac-Moody groups and their flag manifolds.
By extracting the information of cohomology from the Poincar\'{e} series, we are able to determine the rational cohomology rings
of Kac-Moody groups and their flag manifolds. Since Kac-Moody groups and their flag manifolds are rational formal, we also determine their rational homotopy groups and rational homotopy types.
\end{abstract}

\section{Introduction}
Let $A=(a_{ij})$ be an $n\times n$ integer matrix satisfying

(1) For each $i,a_{ii}=2$;

(2) For $i\not=j,a_{ij}\leq 0$;

(3) If $a_{ij}=0$, then $a_{ji}=0$.

\noindent then $A$ is called a Cartan matrix.

By the work of Kac\cite{Kac_68} and Moody\cite{Moody_68} it is well known that for each $n\times n$ Cartan matrix $A$, there is a Lie algebra $g(A)$ associated to
$A$ which is called Kac-Moody Lie algebra. Then Kac and Peterson\cite{Kac_Peterson_83}\cite{Kac_Peterson_84}\cite{Kac_85} constructed the Kac-Moody group $G(A)$ corresponding to the Lie algebra $g(A)$. In this paper we consider
the quotient Lie algebra  of $g(A)$ modulo its center $c(g(A))$ and the associated simply connected group $G(A)$ modulo $C(G(A))$, for convenience we still use the symbols $g(A)$ and $G(A)$. %By the Proposition 1.6 in \cite{Kac_82}, the Cartan subgroup of $G(A)$ has rank $n$.

Kac-Moody Lie algebras and Kac-Moody groups are divided into three types.

(1) Finite type, if $A$ is positive definite. In this case, $G(A)$ is just the simply connected complex semisimple Lie group with Cartan matrix $A$.

(2) Affine type, if $A$ is positive semi-definite and has rank $n-1$.

(3) Indefinite type otherwise.

Denote the Weyl group of ${g}(A)$ by $W(A)$, then $$W(A)=<\sigma_1,\cdots,\sigma_n|\sigma^2_i=1,1\leq  i\leq n; (\sigma_i\sigma_j)^{m_{ij}}=1,1\leq i<j\leq n>.$$ where $\sigma_1,\sigma_2,\cdots,\sigma_n$ are the Weyl reflections with respect to $n$ simple roots $\alpha_1,\alpha_2,\cdots,\alpha_n$ of $g(A)$, $m_{ij}=2,3,4,6 $ and $\infty$ as $a_{ij}a_{ji}=0,1,2,3 $ and $\geq 4$ respectively.

Each element $w\in W(A)$ has a decomposition of the form $w=\sigma_{i_1}\cdots \sigma_{i_k},1\leq i_1,\cdots,i_k\leq n$. The length of $w$ is defined as the least integer $k$ in all of those decompositions of $w$, denoted by $l(w)$. The Poincar\'{e} series of $g(A)$ is the power series $P_A(q)=\sum\limits_{w\in W(A)} q^{2l(w)}$.

For the Kac-Moody Lie algebra $g(A)$, there is the Cartan decomposition $g(A)=h\oplus \sum\limits_{\alpha\in \Delta} g_{\alpha}$, where $h$ is the Cartan sub-algebra and $\Delta$ is the
root system of $g(A)$. Let $b=h\oplus \sum\limits_{\alpha\in \Delta^+} g_{\alpha}$ be the Borel sub-algebra, then $b$ corresponds to a Borel subgroup $B(A)$ in the Kac-Moody group $G(A)$. The
homogeneous space $F(A)=G(A)/B(A)$ is called the flag manifold of $G(A)$. By Kumar\cite{Kumar_02}, $F(A)$ is an ind-variety.

The flag manifold $F(A)$ admits a CW-decomposition of Schubert cells which are indexed by the elements of Weyl group $W(A)$. For each $w\in W(A)$, the real dimension of Schubert variety $X_w$ is
$2l(w)$. So the Poincar\'{e} series of flag manifold $F(A)$ is just the Poincar\'{e} series $P_A(q)$ of $g(A)$.

By the well known results about the cohomology and Poincar\'{e} series of flag manifolds of Kac-Moody groups, for example see Kichiloo\cite{Kitchloo_98} and Kumar\cite{Kumar_02}, it follows that: the rational cohomology rings $H^*(G(A))$ and $H^*(F(A))$ are locally finite. Hence they are generated by countable number of generators. Steinberg\cite{Steinberg_68} proved that
the Poincar\'{e} series $P_A(q)$ of the Kac-Moody flag manifold $F(A)$ is a rational function.

The rational cohomology rings of Kac-Moody groups and their flag manifolds of finite or affine type are extensively studied by many algebraic topologists. For reference, see Pontrjagin\cite{Pontryagin_35}, Hopf\cite{Hopf_41}, Borel\cite{Borel_53_1}\cite{Borel_53}\cite{Borel_54}, Bott and Samelson\cite{Bott_Samelson_55}, Bott\cite{Bott_56}, Milnor and Moore\cite{Milnor_Moore_65}, Chevalley\cite{Chevalley_94}, Bernstein, Gel'fand and Gelfand\cite{BGG_73}, Demazure\cite{D_75}, Kostant and Kumar\cite{Kostant_Kumar_86} and Kichiloo\cite{Kitchloo_98} etc.

%The rational cohomology ring $H^*(G)$ of Lie groups are extensively studied by vast number of algebraic topologists, for example by
%Milnor and Moore\cite{Milnor_Moore_65} etc.

Below is the well known theorem about the rational cohomology of Hopf spaces.

\noindent {\bf Theorem}(Hopf): Let $G$ be a connected H-space which has the homotopy type of a CW-complex, then the rational cohomology rings $H^*(G)$ is a Hopf algebra and as algebra it is isomorphic to the the tensor product of a polynomial algebra $P(V_0)$ and
a exterior algebra $\Lambda(V_1)$, where $V_0$ and $V_1$ are
respectively the set of even and odd degree free generators of $H^*(G)$.

For a Kac-Moody group $G(A)$, $H^k(G(A))$ is a finite dimension rational vector space for each $k\geq 0$. Denote the number of degree $k$ generators in $V=V_0\cup V_1$ by $i_k$, then the Poincar\'{e} series of $G(A)$ is
\begin{equation}\label{1}
 P_G(q)=\prod\limits_{k=1}^{\infty}\frac{ (1-q^{2k-1})^{i_{2k-1}}}{(1-q^{2k})^{i_{2k}}}.
\end{equation}

An important result is

\noindent{\bf Theorem 1: }The cohomology ring $H^*(G(A))$ is determined by its Poincar\'{e} series.%Kac-Moody group $G(A)$ has the homotopy type of a locally finite CW-complexes, and

The Poincar\'{e} series for finite and affine type Kac-Moody groups and their flag manifolds are known.
The computations of Poincar\'{e} series for hyperbolic Kac-Moody flag manifolds are discussed in
Gungormez and Karadayi\cite{GK_arxiv}, Chapovalov, Leites and Stekolshchik\cite{CDR_10}. The general indefinite case is considered by the authors in \cite{Jin_10}\cite{Jin_Zhao_11}.

But for rational cohomology rings of indefinite Kac-Moody groups and their flag manifolds, little is known. In fact except for $n\leq 2$(see \cite{Kitchloo_98}), not a single case is computed. So a natural question is

\noindent{\bf Problem: }Compute the rational cohomology rings for Kac-Moody groups and their flag manifolds of indefinite type.

In this paper our main goal is to solve this problem. Since Kac-Moody groups and their flag manifolds are rational formal\cite{Sullivan_77}, their rational homotopy types are determined by their rational cohomology rings. Our strategy is:

(1) Compute the Poincar\'{e} series $P_A(q)$ of $F(A)$ from the Poincar\'{e} series of $G(A)$ by applying Leray-Serre spectral sequence;

(2) Compute the Poincar\'{e} series $P_A(q)$ of $F(A)$ from the Cartan matrix $A$ by using the method in \cite{Jin_10}\cite{Jin_Zhao_11};

(3) Determine $i_k,k\geq 1$ in Equation (1) by comparing
the Poincar\'{e} series obtained in (1) and (2).

%then using the Leray-Serre spectral sequence to the fibre bundle
%$G(A)\to F(A)\to BB(A)$ to compute the Poincar\'{e} series of $G$. analysting

Let the classifying map of principal $B(A)$-bundle
$\pi:G(A)\to F(A)$ be $j: F(A)\to BB(A)$, where $BB(A)$ denotes the classifying space of Borel subgroup $B(A)$. There is a homotopy fibration $G(A)\stackrel{\pi}{\longrightarrow} F(A) \stackrel{j}{\longrightarrow} BB(A)$. By analyzing the Leray-Serre spectral sequence, we show

\noindent {\bf Theorem 2: }For a Kac-Moody group $G(A)$ with Poincar\'{e} series $P_G(q)$ as in equation (1), the Poincar\'{e} series of $F(A)$ is
\begin{equation}\label{2}
       P_A(q)=\frac{\prod\limits_{k=1}^{\infty} (1-q^{2k})^{i_{2k-1}}}{(1-q^2)^n}   \frac{1}{\prod\limits_{k=1}^{\infty}(1-q^{2k})^{i_{2k}}}.
\end{equation}
%\noindent here $n$ is the rank of $G(A)$.

For a Kac-Moody group $G(A)$, it is shown by Kac\cite{Kac_852} that the set $V_1$ is a finite set and the number $l$ of elements in $V_1$ is less or equal to $n$. So the exterior algebra part of $H^*(G(A))$ is of finite dimension.

%Let $m=n-l$, then we write the Poincar\'{e} series of $F(A)$ into product of two parts. Let $$Q_A(q)= \frac{\prod\limits_{k=1}^{\infty} (1-q^{2k+2})^{i_{2k+1}+1}}{(1-q^2)^{l}}$$
%and $$R_A(q)={(1-q^2)^m }{\prod\limits_{k=1}^{\infty}(1-q^{2k})^{i_{2k}}},$$
%then $P_A(q)=Q_A(q)/R_A(q)$. $R_A(q)$ corresponds to
%the polynomial part of the rational cohomology of $F(A)$ and $Q_A(q)$ corresponds to
%the finite part. A surprising result is

%{\bf Theorem 3: } In fact $R_A(q)$ is a polynomial of $q$.

%In fact $Q_A(q)$ cothe Poincar\'{e}  product of the finite part infinite part.

\noindent{\bf Theorem 3:} Let $P_A(q)$ be the Poincar\'{e} series of flag manifold $F(A)$, then the sequence $i_2-i_1,i_4-i_3,\cdots,i_{2k}-i_{2k-1},\cdots$ can be derived from $P_A(q)$. In fact $P_A(q)$ can also be recovered from the sequence $i_2-i_1,i_4-i_3,\cdots,i_{2k}-i_{2k-1},\cdots$.

But to determine the rational homotopy type and rational cohomology of $G(A)$, we need to determine the sequence $i_1,i_2,\cdot\cdot\cdot,i_k,\cdot\cdot\cdot$. So except for the Poincar\'{e} series $P_A(q)$, more ingredients are needed. If we can determine all the degrees of elements in $V_1$, then we will determine the sequence $i_1,i_3,\cdot\cdot\cdot,i_{2k-1},\cdot\cdot\cdot$. Combining with Theorem 3, we also determine the sequence $i_2,i_4,\cdot\cdot\cdot,i_{2k},\cdot\cdot\cdot$, hence work out the rational cohomology ring and rational homotopy type of $G(A)$.

%So if we can further determine all the degrees of elements in $V_1$, then we can determine the $i_1,i_3,\cdots,i_{2k-1},\cdots $. Combining with Theorem 3, we also determine $i_2,i_4,\cdots,i_{2k},\cdots $, hence determine the rational cohomology of $G(A)$.

By Kac\cite{Kac_852}, Kac and Peterson\cite{Kac_Peterson_853} the sequence $i_1,i_3,\cdots,i_{2k-1},\cdots $ can be determined by the ring of polynomial invariants of Weyl group $W(A)$. In fact $i_{2k-1}$ is just the number of degree $k$ basic invariant polynomials of $W(A)$.

In \cite{Zhao_Jin_12} the authors prove the following result.

\noindent {\bf Theorem 4: }Let $A$ be an indecomposable and indefinite Cartan matrix. If $A$ is symmetrizable, then $I(A)=\mathbb{Q}[\psi]$; If $A$ is non symmetrizable, then $I(A)=\mathbb{Q}$. Where $I(A)$ is the ring of polynomial invariants of Weyl group $W(A)$, and $\psi$ is a invariant bilinear form on $g(A)$.

This theorem is a generalization of Moody's result(see \cite{Moody_78}) for symmetrizable hyperbolic Cartan matrices.

By Theorem 4, for an indecomposable and indefinite Cartan matrix $A$, $i_{2k-1}=0$ for all $k>0$ except for $k=2$. And for $k=2$, if $A$ is symmetrizable, $i_3=1$; If $A$ is non symmetrizable, $i_3=0$. Hence the sequence $i_1,i_3,\cdots,i_{2k-1},\cdots $ is worked out. Combining with Theorem 3 we eventually determine the sequence $i_1,i_2,\cdots, i_{k},\cdots$. As a consequence the rational cohomology and rational homotopy types of the Kac-Moody group $G(A)$ and its flag manifold $F(A)$ are determined.

This paper is organized as below. In section 2, we introduce some algebraic and combinatorial results we needed. %show for any polynomial $f(q)$ with integer coefficients and constant item 1, how it be rewrited into the form of product $\prod\limits_{k=1}^{\infty}(1-q^k)^{i_k}$.
In section 3 we discuss the decompositions of order of Weyl group $W(A)$ and Poincar\'{e} series of flag manifold $F(A)$. In section 4 we compute the Poincar\'{e} series $P_A(q)$ of flag manifold $F(A)$ from the Poincar\'{e} series $P_G(q)$ of $G(A)$ via Leray-Serre spectral sequence. The main theorem of this paper will be concluded in Section 5. And an example is given in section 6.

\section{Some algebras and combinatorics}

Let $\mathbb{Z}_1[q]$ be the set of power series of $q$ with integer coefficients and the constant item $1$. That is for $f(q)\in \mathbb{Z}_1[q]$,
$$f(q)=1+a_1 q+a_2q^2+\cdots +a_k q^k+\cdots,a_k\in \mathbb{Z},k>0.$$

We have the following result.

\noindent {\bf Proposition 1: }$f(q)\in \mathbb{Z}_1[q]$ can be expanded uniquely into the form of products ${\prod\limits_{k=1}^{\infty}(1-q^{k})^{i_{k}}}$, where
$i_1,i_2,\cdots,i_k,\cdots$ are integer sequence.

\noindent {\bf Proof:} For $f(q)$, we define an integer sequence $i_1,i_2,\cdots,i_k,\cdots$ inductively.

At first we set $i_1=a_1$. Let $f^{(1)}(q)=f(q)/(1-q)^{i_1}$, then $f^{(1)}\in \mathbb{Z}_1[q]$ and
$$f^{(1)}(q)=1+ a_2^{(1)}q^2+a_3^{(1)}q^3+\cdots +a_k^{(1)} q^k+\cdots .$$
Set $i_2=a_2^{(1)}$ and let $f^{(2)}(q)=f^{(1)}(q)/(1-q^2)^{i_2}$, then $f^{(2)}\in \mathbb{Z}_1[q]$ and
$$f^{(2)}(q)=1+ a_3^{(2)}q^3+a_4^{(2)}q^4+\cdots +a_k^{(2)} q^k+\cdots .$$
Set $i_3=a_3^{(2)}$ and continue the same procedure we get integer sequence $i_1,i_2,\cdots,i_k,\cdots$.

By checking the previous procedure and comparing the coefficients of the power series concerned, we prove the uniqueness of the sequence $i_1,i_2,\cdots,i_k,\cdots$.

We call the sequence $i_1,i_2,\cdots,i_k,\cdots$ the characteristic sequence of power series $f(q)$. By Proposition 1, there is a one to one
correspondence between $\mathbb{Z}_1[q]$ and the set of integer sequences indexed on the set $\mathbb{N}$ of natural numbers. And if $f_1(q),f_2(q)$ correspond to the integer sequences $i_1,i_2,\cdots,i_k,\cdots$ and $j_1,j_2,\cdots,j_k,\cdots$, then the product $f_1(q)f_2(q)$ corresponds to
integer sequence $i_1+j_1,i_2+j_2,\cdots,i_k+j_k,\cdots$.

The Poincar\'{e} series $P_A(q)$ of $F(A)$ is of the form $P_A(q)=G(q^2)$. By applying Proposition 1 to $G(q)$ we will prove Theorem 1.

%\noindent {\bf Corollary 1: } $f(q)\in \mathbb{Z}_1[q]$ is the Poincar\'{e} series of a connected locally finite free graded algebra over rational $\mathbb{Q}$ generated by polynomial and %exterior algebra generators if and only if its corresponding integer sequence take values in $\mathbb{Z}^+=\{n\in \mathbb{Z}|n\geq 0\}$.

The characteristic sequence $i_1,i_2,\cdots,i_k,\cdots$ can be computed by using M\"{o}bius Inversion Theorem.

\noindent {\bf Proposition 2: } Let $f(q)\in \mathbb{Z}_1[q]$, suppose $\ln( f(q))=b_1 q+\frac{b_2}{2} q^2 \cdots+\frac{b_k}{k} q^k+\cdots$, then $$i_k=\frac{1}{k}\sum\limits_{n|k}\mu(n)b_{k/n}.$$
\noindent {\bf Proof: }For $$f(q)=\prod\limits_{k=1}^{\infty}(1-q^k)^{i_k}$$
take \lq\lq$\ln$" to the two sides of equation $\displaystyle{f(q)=\prod\limits_{k=1}^{\infty}(1-q^k)^{i_k}}$, we get
$$b_1 q+\frac{b_2}{2} q^2 \cdots+\frac{b_k}{k} q^k+\cdots=\ln\prod\limits_{k=1}^{\infty}(1-q^k)^{i_k}$$
the right side$=\displaystyle{\sum\limits_{k=1}^{\infty}i_{k}\ln(1-q^k)
   =\sum\limits_{k=1}^{\infty}i_k\sum\limits_{n=1}^{\infty}\frac{q^{kn}}{n}=\sum\limits_{k=1}^{\infty}\sum\limits_{n=1}^{\infty}\frac{i_k}{n}q^{kn}=\sum\limits_{n=1}^{\infty}(\sum\limits_{k|n}
   \frac{k\cdot i_k}{n})q^n }$

By comparing the coefficients of the two sides we get $\sum\limits_{k|n}k\cdot i_k=b_n$

We need the following theorem to compute $i_k$.

\noindent\textbf{ M\"{o}bius Inversion Theorem:} Let $F(n),f(n)$ be two sequences indexed on natural numbers, if $F(n)=\sum\limits_{k|n}f(k)$, then $f(k)=\sum\limits_{n|k}\mu(n)F(k/n)$, where $\mu(n)$ satisfies:

(1) When $n=1,\mu(n)=1$.

(2) When $n=\prod\limits_{i=1}^r p_i$, $\mu(n)=(-1)^r$, where $p_i$'s are different prime integers.

(3) $\mu(n)=0$ otherwise.

By the M\"{o}bius Inversion Theorem, from $\sum\limits_{k|n}k\cdot i_k=b_n$ we get
$$k\cdot i_k=\sum\limits_{n|k}\mu(n)b_{k/n}$$
\noindent {\bf Example 1: }If $f(q)=1-2q$, then
in this case, $b_n=2^n $, so
$$k\cdot i_k=\sum\limits_{d|k}\mu(d)2^{\frac{k}{d}}.$$
Take $k=18$, all the divisors of $18$ is $1,2,3,6,9,18$ and $$\mu(1)=1,\mu(2)=-1,\mu(3)=-1,\mu(6)=1,\mu(9)=\mu(18)=0,$$
then $18i_{18}=2^{18}-2^9-2^6+2^3$, so $i_{18}=14532$.

We list the computation results as below for $1\leq k\leq 18$.
$$\begin{tabular}{|c|c|c|c|c|c|c|c|c|c|c|c|c|c|c|c|c|c|c|}
  \hline
  $k$ & 1 & 2 & 3 & 4 & 5 & 6 & 7 & 8 & 9 & 10 & 11 & 12 & 13 & 14 & 15 & 16 & 17 & 18  \\ \hline
  $i_k$ & 2 & 1 & 2 & 3 & 6 & 9 & 18 & 30 & 56 & 99 & 186 & 335 & 630 & 1161 & 2182 & 4080 & 7710 & 14532  \\
  \hline
\end{tabular}$$

\section{Decompositions of Weyl groups and Poincar\'{e} series}
Let $G(A)$ be a Kac-Moody group and $W(A)$ be its Weyl group. From the definition of Poincar\'{e} series $P_A(q)$, we see as $q$ approaches to $1$, $P_A(q)$ tends to $|W(A)|$(the order of $W(A)$). Hence the Poincar\'{e} series $P_A(q)=\sum\limits_{w\in W(A)} q^{2l(w)}$ can be regarded as deformation(or quantization) of $|W(A)|$.

For classical Lie group of type $A_n$, $G(A_n)=\mathrm{SL}(n+1,\mathbb{C})$. The Weyl group of $G(A_n)$ is $W(A_n)\cong S_{n+1}$, the permutation group of set $\{1,2,\cdots,n,n+1\}$,
$|W(A_{n})|=(n+1)!$. The Poincar\'{e} series of flag manifold $F(A_n)=\mathrm{SL}(n+1,\mathbb{C})/B$ is $P_{A_n}(q)=\prod\limits_{i=1}^{n+1} \displaystyle{\frac {1-q^{2k}}{1-q^2}}$. If we set $[k]=\displaystyle{\frac{1-q^{2k}}{1-q^2}}$, then as $q$ approaches
to $1$, $[k]$ approaches to $k$. The Poincar\'{e} series is $[n+1]!=[n+1][n]\cdots[1]$.

In fact for the Kac-Moody group $G(A)$ of finite type, if $$H^*(G(A))\cong \Lambda(x_1,x_2,\cdots,x_n), \deg x_k=2d_k-1,$$
then $|W(A)|=\prod\limits_{k=1}^n d_k$ and $P_A(q)=\prod\limits_{k=1}^n [d_k]$, where $d_k$'s are the degrees of the basic invariant polynomials of $W(A)$.

%The order of the Weyl group $|W(A)|$ of a Kac-Moody group $G(A)$ of finite type is finite.
For a Kac-Moody group $G(A)$ of affine or indefinite type, $|W(A)|$ is infinite, so at first sight the decomposition of $|W(A)|$ is meaningless.
But if we consider $P_A(q)$, %the quantization of $|W(A)$,
through a regularization procedure as in quantum fields theory, we can give an interesting decomposition.

%each Kac-Moody group $G(A)$ of finite type, there is
Let $G(\widetilde A)$ be an untwisted affine Kac-Moody groups with $\widetilde A$ the extended Cartan
matrix of a Cartan matrix $A$ of finite type. For the flag manifold $F(\widetilde A)$, the Poincar\'{e} series is
$$P_{\widetilde{A}}(q)=P_A(q)\prod\limits_{i=1}^{n}\frac{1}{1-q^{2d_k-2}}=\prod\limits_{k=1}^n [d_k] [\infty]_{d_k-1}.$$
where $[\infty]_k=\frac {1}{1-q^{2k}}$.

So we have the decomposition $|W(\widetilde A)|=\prod\limits_{k=1}^n d_k \prod\limits_{k=1}^n\infty_{d_k-1}$.

%For untwisted affine flag manifolds we separated the finite and infinite part of $|W(A)|$ successfully.
For the Weyl groups and Poincar\'{e} series of twisted affine flag manifolds, the computation results are almost the same as the twisted case.
The decomposition of order of Weyl groups and Poincar\'{e} series for indefinite case can also be done. But due to the relation $[k]=[\infty]_1 [\infty]_k^{-1}$, the decomposition is not unique. For details see \cite{Jin_Zhao_11}.

%In \cite{Jin_Zhao_11} we consider the computation of Poincar\'{e} series $P_A(q)$. The computation results shows the Poincar\'{e} series split to two parts,
%the finite part and infinite part.

%Unlike the case for Poincar\'{e} series of $G(A)$, the Poincar\'{e} series $P_A(q)$ of flag manifold $F(A)$ can't determine
%$$i_1,i_2,\cdots,i_k,\cdots$ uniquely. Suppose $$
% P_A(q)=[\infty]_1^n\frac{\prod\limits_{k=1}^{\infty} (1-q^{2k})^{i_{2k-1}}}{(1-q^2)^n}   \frac{1}{\prod\limits_{k=1}^{\infty}(1-q^{2k})^{i_{2k}}}
%$$that is $$P_A(q)=\prod\limits_{k=1}^{\infty} [k]^{i_{2k-1}}\prod\limits_{k=1}^{\infty}[\infty]_k^{i_{2k}}.$$
%simplify this equation with the formula $[k]=[\infty]_1 [\infty]_k^{-1}$, we get
%\begin{equation}%\label{}
% P_A(q)=\frac{1}{(1-q^2)^{n+i_2-i_1}}   \frac{1}{\prod\limits_{k=2}^{\infty}(1-q^{2k})^{i_{2k}-i_{2k-1}}}
%\end{equation}

%By using Proposition 1 to equation (3), we have

%\noindent {\bf Theorem 3: }The sequence $n+i_2-i_1,i_4-i_3,\cdots,i_{2k}-i_{2k-1},\cdots$ is uniquely determined by Poincar\'{e} series $P_A(q)$.

%It is obvious that the sequence $i_2-i_1,i_4-i_3,\cdots,i_{2k}-i_{2k-1},\cdots$ also determines $P_A(q)$, so it contains all the information we can get from $P_A(q)$.

\section{Leray-Serre spectral sequence}
In this section we use the Leray-Serre spectral sequence of fibration $G(A)\stackrel{\pi}{\to} F(A)\stackrel{j}{\to} BB(A)$ to compute the cohomology and Poincar\'{e} series of $F(A)$ from those of $G(A)$.
For reference see \cite{Whitehead_78}.

$BB(A)$ is homotopic to the $n$-fold Cartesian product of classifying space of $B \mathbb{C}^*$($\mathbb{C}^*=\mathbb{C}-\{0\}$), denote the cohomology generators of $H^*(BB(A))$ by $\omega_1,\cdots,\omega_n,\deg \omega_i=2$, the free generators of $V_1$ by $y_1,\cdots,y_l$, and the free generators of $V_0$ by $z_1,\cdots,z_k,\cdots$. Where $\omega_1,\cdots,\omega_n$ correspond to the fundamental dominant weights of $g(A)$. We have the spectral sequence $(E_r^{p,q},d_r)$ with
$$E_2^{p,q}=H^p(BB(A);H^q(G(A)))
%\cong H^p(BB(A))\otimes H^q(G(A))
\cong \mathbb{Q}[\omega_1,\cdots,\omega_n]\otimes\Lambda(y_1,\cdots,y_l)\otimes \mathbb{Q}[z_1,\cdots,z_k,\cdots].$$
The differential $d_2: E_2^{p,q}\to E_2^{p+2,q-1}$ is given by $d_2(\omega_i\otimes 1)=0,d_2(1\otimes z_k)=0$ and $d_2(1\otimes y_j)=f_j\otimes 1, 1\leq j\leq l$, where $f_j\in H^*(BB(A))$ is a polynomial of $\omega_i,1\leq i\leq n$. %Hence by the property of the differential $d_2(t_i\otimes y_j)=t_if_j\otimes 1$.
A routine computation shows $$H^*(F(A))\cong E_3^{*,*}\cong \mathbb{Q}[\omega_1,\cdots,\omega_n]/<f_j,1\leq j\leq l>\otimes \mathbb{Q}[z_1,\cdots,z_k,\cdots].$$

By Kac\cite{Kac_852}, $f_1,\cdots,f_{l}$ is a regular sequence in $H^*(BB(A))$, so $l\leq n$. Therefore we conclude that the Poincar\'{e} series of $F(A)$ is given by the Formula (2). This proves the Theorem 2.

By simplifying the Formula (2), we get
\begin{equation}%\label{}
 P_A(q)=\frac{1}{(1-q^2)^{n+i_2-i_1}}   \frac{1}{\prod\limits_{k=2}^{\infty}(1-q^{2k})^{i_{2k}-i_{2k-1}}}
\end{equation}

By using Proposition 1 to Equation (3), we prove the Theorem 3.

\section{The rational homotopy type of $G(A)$}

The Poincar\'{e} series $P_A(q)$ can be computed easily by an inductive procedure, see \cite{Jin_10}\cite{Jin_Zhao_11} for example.
From $P_A(q)$ we can compute the sequence $i_2-i_1,i_4-i_3,\cdots,i_{2k}-i_{2k-1},\cdots$. But to determine the rational homotopy type of $G(A)$, we need to determine the sequence $i_1,i_2,\cdots,i_k,\cdots$. So except for the Poincar\'{e} series $P_A(q)$, we need more ingredients. Theorem 4 fills the gap. For its proof, see \cite{Zhao_Jin_12}.

By Theorem 4 for indecomposable and indefinite Cartan matrix $A$, $i_1,i_3,\cdots,i_{2k+1},\cdots $ is determined, combining with Theorem 3, we determine the complete sequence $i_1,i_2,\cdots,i_k,\cdots$.

Since $G(A)$ is a simply connected group, we get $i_1=0$. By virtue of Schubert decomposition, we know that $H^2(F(A))$ is spanned by the Schubert classes corresponding to $n$ Weyl
reflections. So $n+i_2-i_1=n$, hence $i_2=0$.
To determine $i_3$ we need the following definition which is well known in Kac-Moody Lie algebras theory.

\noindent {\bf Definition 1:} An $n\times n$ Cartan matrix $A$ is symmetrizable if there exist an invertible diagonal matrix $D$ and a symmetric matrix $B$
such that $A=DB$. $g(A)$ is called a symmetrizable Kac-Moody Lie algebra if its Cartan matrix is symmetrizable.

The symmetrizability of a Cartan matrix $A$ is intimately related to the existence of non degenerate invariant bilinear form $\psi$ on $g(A)$. From Theorem 4, it follows directly that If $g(A)$ is a symmetrizable Kac-Moody Lie algebra, then $i_3=1$, otherwise $i_3=0$. And $i_{2k-1}=0$ for $k\geq 3$.

Set $\epsilon(A)=1$ or $0$ depending on $A$ is symmetrizable or not as in \cite{Kac_852}, then we can state our main theorem as below.

%\noindent {\bf Theorem:} (Moody) Let A be an indecomposable symmetrizable $n\times n$
%Cartan matrix whose associated invariant quadratic form $\psi$ is non-degenerate and of signature
%$(n-1,1)$, then the ring of $W(A)$ invariant polynomial functions on $\mathbb{Q}$ is $\mathbb{Q}[\psi]$.

%In his paper "Polynomial invariants of isometry groups of indefinite quadratic lattices", Moody further said:

%{\sl It is also true for a wide class of Weyl
%groups of infinite root systems, including all the hyperbolic root systems
%and we conjecture that it is in fact true for all
%Weyl groups arising from non-singular Cartan matrices of non-finite
%type.}

%In his paper "Constructing groups associated to infinite-dimensional Lie algebras", Kac also use the result that the assumption that for indefinite and indecomposable Cartan matrix
%the ring of $W(A)$ invariant polynomial functions on $\mathbb{Q}$ is $\mathbb{Q}[\psi]$.

%Though it is not clear in what extent "The ring of $W(A)$ invariant polynomial functions on $\mathbb{Q}$ is $\mathbb{Q}[\psi]$" is true, but we can believe that it is valid for a wide class of
%indefinite Cartan matrices. So we called the property "The ring of $W(A)$ invariant polynomial functions on $\mathbb{Q}$ is $\mathbb{Q}[\psi]$" by Property V.

\noindent{\bf Theorem 5: }For an indecomposable and indefinite Cartan matrix $A$, if $A$ is symmetrizable, then
$$H^*(G(A))\cong \Lambda_{\mathbb{Q}}(y_3)\otimes \mathbb{Q}[z_1,\cdots,z_k,\cdots]$$ and $$H^*(F(A))\cong  \mathbb{Q}[\omega_1,\cdots,\omega_n]/<\psi>\otimes \mathbb{Q}[z_1,\cdots,z_k,\cdots]. $$
If $A$ is non symmetrizable, then
$$H^*(G(A))\cong \mathbb{Q}[z_1,\cdots,z_k,\cdots]$$ and $$H^*(F(A))\cong  \mathbb{Q}[\omega_1,\cdots,\omega_n]\otimes \mathbb{Q}[z_1,\cdots,z_k,\cdots]. $$
where $\deg z_k\geq 4$ are even for all $k$ and their degrees can be determined from the Poincar\'{e} series $P_A(q)$ and $\epsilon(A)$.

%{\bf Theorem 4: }For a Cartan matrix $A$ satisfies Property V, the sequence $i_1,i_2,\cdots,i_k,\cdots$ is determined by the Poincar\'{e} series $P_A(q)$ and the symmetrizability of $A$. The
%cohomology ring $$H^*(G(A))\cong \Lambda_{\mathbb{Q}}(y_3)\otimes \mathbb{Q}[z_1,\cdots,z_k,\cdots]$$ and $$H^*(F(A))\cong  \mathbb{Q}[t_1,\cdots,t_n]/<\psi>\otimes %\mathbb{Q}[z_1,\cdots,z_k,\cdots]. $$

%Therefore for a Cartan matrix $A$ such that Property V is satisfied,
The rational homotopy groups and the rational minimal model of the Kac-Moody group $G(A)$ and its flag manifold $F(A)$ can be computed from this theorem easily.
%e Cartan matrix $A$. And the rational homotopy of $G(A)$ and $F(A)$ is determined.

%{\bf Corollary: }

%cohomology

%the minimal model

%the rational homotopy groups
Kumar\cite{Kumar_85} proved that for Kac-Moody Lie algebra $g(A)$, the Lie algebra cohomology $H^*(g(A),\mathbb{C})$ $\cong H^*(G(A))\otimes \mathbb{C}$, so we also compute $H^*(g(A),\mathbb{C})$.

\section{An example}
We end this paper by the following example.

\noindent{\bf Example 2: }Let $A$ be a rank $n$ Cartan matrix $A$ which satisfies $a_{ij}a_{ji}\geq 4$ for all $i\not =j$, then the Weyl group of $G(A)$ is
$$W(A)=<\sigma_1,\cdots,\sigma_n|\sigma^2_i=1,1\leq i\leq n>.$$
The length $k$ elements in $W(A)$ are in one to one correspondence to the words $\sigma_{i_1}\sigma_{i_2}\cdots \sigma_{i_k}$ satisfying $i_t\not= i_{t+1}, 1\leq t\leq k-1$. So the number of length $k$ elements in $W(A)$ is $n(n-1)^{k-1}$ and the Poincar\'{e} series is $$P_A(q)=\sum\limits_{k=0}^{\infty} n(n-1)^{k-1}q^{2k}=\frac{1+q^2}{1-(n-1)q^2}$$

By the well known Witt formula(see \cite{Kang_93}) $$1-nq=\prod\limits_{k=1}^\infty (1-q^{k})^{\dim L^k_n}$$ where $L_n$ is the free graded Lie algebra generated by $n$ elements with degree $1$ and $L_n=\bigoplus\limits_{k=1}^\infty L_n^k$. $\dim L_n^k$ is the dimension of the degree $k$ homogeneous component of $L_n$. It is computed by using the M\"{o}bius Inversion Theorem. See Example 1 for the case $n=2$.
Hence $$P_A(q)=\frac{1-q^4}{(1-q^{2})^n}\frac{1}{\prod\limits_{k=2}^\infty (1-q^{2k})^{\dim L^k_{n-1}}}.$$

From the expression of $P_A(q)$, we get:

If $A$ is symmetrizable, then $i_1=i_2=0,i_3=1$ and $i_{2k-1}=0$ for $k\geq 3$, $i_{2k}=\dim L^k_{n-1}$ for $k\geq 2$.

If $A$ is not symmetrizable, then $i_1=i_2=i_3=0, i_4=\dim L^2_{n-1}-1$ and $i_{2k-1}=0 ,i_{2k}=\dim L^k_{n-1}$ for $k\geq 3$.

\end{document}